\documentclass[reqno]{amsart}
\usepackage{a4wide}
\usepackage{todonotes}
\usepackage[all,2cell]{xy}
\usepackage{graphicx, amsmath, amssymb, amsthm}
\usepackage{newclude}
\usepackage{hyperref}
\usepackage{enumitem}
\usepackage[utf8]{inputenc}
\usepackage{staves}

\newcommand{\<}{\langle}
\let\>=\undefined
\newcommand{\>}{\rangle}

\newcommand{\st}{\mid}

\newcommand{\ds}{\displaystyle}

\newcommand{\bdry}{\partial}

\newcommand{\tnsr}{\otimes}

\newcommand{\morph}{\mathop{\longrightarrow}\limits}
\newcommand{\comorph}{\mathop{\longleftarrow}\limits}
\newcommand{\twomorph}{\mathop{\Longrightarrow}\limits}

\newcommand{\from}{\colon}
\newcommand*\lon{%
        \nobreak
        \mskip6mu plus1mu
        \mathpunct{}%
        \nonscript
        \mkern-\thinmuskip
        {:}%
        \mskip2mu
        \relax
}

\renewcommand{\projlim}{\mathop{\operatorname*{lim}\limits_{\longleftarrow}}\limits}
\newcommand{\colim}{\mathop{\operatorname*{lim}\limits_{\longrightarrow}}\limits}

\newcommand{\mono}{\hookrightarrow}
\newcommand{\monic}{\mathop{\rightarrowtail}\limits}

\newcommand{\isom}{\cong}
\newcommand{\adjoint}{\dashv}

\newcommand{\cA}{\mathcal A}

\newcommand{\cC}{\mathcal C}
\newcommand{\cE}{\mathcal E}

\newcommand{\cH}{\mathcal H}

\newcommand{\cO}{\mathcal O}

\newcommand{\cV}{\mathcal V}

\newcommand{\cX}{\mathcal X}

\newcommand{\frC}{\mathfrak C}
\newcommand{\frD}{\mathfrak D}

\newcommand{\frV}{\mathfrak V}
\newcommand{\frW}{\mathfrak W}

\newcommand{\bbK}{\mathbb K}

\newcommand{\bbQ}{\mathbb Q}
\newcommand{\bbR}{\mathbb R}

\newcommand{\bbT}{\mathbb T}

\newcommand{\mlaux}[3]{\setbox0=\hbox{$\mathsurround=0pt #2{#3}$}%
  \dimen0=\dp0\advance\dimen0 by \ht0\lower#1\dimen0\box0}
\newcommand{\mlower}[2]{\mathpalette{\mlaux{#1}}{#2}}
\newcommand{\longtwoheadrightarrow}{\mathrel{\mathord{-}\mkern-3mu\mathord\twoheadrightarrow}}

\renewcommand{\hom}{\underline{\operatorname{hom}}}
\newcommand{\map}{\underline{\operatorname{map}}}
\newcommand{\fibr}{\mathop{\twoheadrightarrow}\limits}
\newcommand{\acyc}{\overset{\mlower{-1}\sim}{\smash\longrightarrow}}
\newcommand{\afibr}{\overset{\mlower{-1}{\sim}}{\smash\longtwoheadrightarrow}}

\newcommand{\slot}{{\blacksquare}}

\newcommand{\cf}{{\mathrm{cf}}}

\newcommand{\op}{\mathrm{op}}
\newcommand{\oop}{\mathrm{du}}

\newcommand{\id}{\mathrm{id}}

\newcommand{\dom}{\mathrm{dom}}
\newcommand{\cod}{\mathrm{cod}}

\newcommand{\lcot}{\mathbin{\hat\pitchfork}}

\DeclareMathOperator{\lan}{\mathsf{lan}}

\DeclareMathOperator{\res}{\mathsf{res}}

\DeclareMathOperator{\Map}{Map}

\newcommand{\term}{1}

\newcommand{\comma}[2]{({#1}\mathbin\downarrow{#2})}

\newcommand{\wLim}[3]{\{#1, #2\}_{#3}}

\newcommand{\Adjf}{{\underline f}}
\newcommand{\Adju}{{\underline u}}
\newcommand{\Adjt}{{\underline t}}

\newcommand{\defeq}{\mathrel{:=}}

\DeclareMathOperator{\ev}{ev}

\newcommand{\resol}[1][g]{{
  {#1}_\bullet
}}
\newcommand{\tresol}{{g^t_\bullet}}

\newcommand{\weirdResol}{{\tilde\ell}}
\newcommand{\magicAdj}{{\ell}}

\newcommand{\coresol}[1][t]{{
  {#1}^\bullet
}}

\newcommand{\compl}[2][t]{{
  {#2}^\wedge_{#1}
}}
\newcommand{\cocompl}[2][g]{{
  {#2}^\vee_{#1}
}}

\newcommand{\bimaps}[5]{\mathchoice{
  \xymatrix@1{
    #2 \ar@<1ex>[r]^-{#4} \ar@{}[r]|-{#1} & #3 \ar@<1ex>[l]^-{#5}
  }
}{
  \xymatrix@1{
    \rule[-1ex]{0pt}{3.5ex}
    #4 \from #2 \ar@<1ex>[r] \ar@{}[r]|-{#1} & #3 \lon #5 \ar@<1ex>[l]
  }
}{}{}
}

\newcommand{\adjnctn}[4]{\bimaps\bot{#1}{#2}{#3}{#4}}

\newcommand{\dadjnctn}[6][2em]{
  \xymatrix@1@C=#1{
    #2 \ar@<2ex>[r]^-{#4} \ar@<-2ex>[r]_-{#6}
    \ar@{}@<1.2ex>[r]|-{\scriptscriptstyle\bot} \ar@{}@<-1.2ex>[r]|-{\scriptscriptstyle\bot}
    &
    #3 \ar[l]|-{#5}
  }
}

\newcommand{\dimorph}[4]{\xymatrix@1{
  #1 \ar@<1ex>[r]^-{#3} \ar@<-1ex>[r]_-{#4} & #2
}}
\newcommand{\twoCell}[5]{\xymatrix@1{
  #1 \ar@<1ex>[r]^-{#3} \ar@{}[r]|{\Downarrow #5} \ar@<-1ex>[r]_-{#4} & #2
}}

\newcommand{\ndo}[2]{\xymatrix@1{
  #1 \ar[r]^-{#2} & #1
}}

\newcommand{\catname}[1]{\mathbf{#1}}
\newcommand{\Adj}{\catname{Adj}}
\newcommand{\Cat}{\catname{Cat}}
\newcommand{\Mnd}{\catname{Mnd}}
\newcommand{\Cmd}{\catname{Cmd}}

\newcommand{\sSet}{\catname{sSet}}
\newcommand{\qCat}{\catname{qCat}}

\newcommand{\Alg}[2]{{#1}^{#2}}
\newcommand{\CoAlg}[2]{{#1}_{#2}}
\newcommand{\Desc}[2]{\frD^{#2}_{#1}}
\newcommand{\CoDesc}[2]{\frC^{#1}_{#2}}

\newcommand{\spc}{{\Delta^\op}}
\newcommand{\aspc}{{\Delta^\op_+}}

\newcommand{\htc}[1]{\catname{Ho}({#1})}

\mathchardef\mhyphen=45

\newtheorem{theorem}{Theorem}[section]
\newtheorem*{theorem*}{Theorem}

\newtheorem{lemma}[theorem]{Lemma}
\newtheorem{corollary}[theorem]{Corollary}
\newtheorem{proposition}[theorem]{Proposition}

\theoremstyle{definition}
\newtheorem{definition}[theorem]{Definition}
\newtheorem{remark}[theorem]{Remark}
\newtheorem{litcomp}[theorem]{Parallel}
\newtheorem{notation}[theorem]{Notation}
\newtheorem{example}[theorem]{Example}
\newtheorem{warning}[theorem]{Warning}

\newcommand{\Ho}{\textrm{Ho}}

\newdir{ >}{{}*!/-10pt/@{>}}
\newcommand{\pullback}{\ar@{}[dr]|<<{\mbox{\Huge$\lrcorner$}}}
\newcommand{\pushout}{\ar@<2pt>@{}[ul]|<{\mbox{\Huge$\ulcorner$}}}

\makeatletter
\providecommand\@dotsep{5}
\renewcommand{\listoftodos}[1][\@todonotes@todolistname]{%
  \@starttoc{tdo}{#1}}
\makeatother

\SelectTips{cm}{}
\UseAllTwocells
\NoCompileMatrices

\begin{document}

\title{$\infty$-categorical monadicity and descent}
\author[Y.~J.~F.~Sulyma]{Yuri~J.~F. Sulyma}
\address{University of Texas \\ Austin, TX 78712}
\email{ysulyma@math.utexas.edu}
\thanks{Sulyma was supported in part by NSF grant DMS-1564289}

\begin{abstract}
Riehl and Verity have introduced an ``$\infty$-cosmic'' framework in which they redevelop the category theory of $\infty$-categories using 2-categorical arguments. In this paper, we begin with a self-contained review of the parts of their theory needed to discuss adjunctions and monadicity. This is applied in order to extend to the $\infty$-categorical context the classical criterion for fully faithfulness of the comparison functor induced by an adjunction. We discuss the relation with previous work in the literature---which primarily uses model-categorical techniques---and indicate applications to descent theory.
\end{abstract}

\maketitle
\tableofcontents


\section{Introduction}

Descent theory plays an important role in algebraic geometry, as well as in the plethora of fields which draw upon its technology. Motivated by the problem of assembling local data into global data, it is profitably reinterpreted in terms of co/monads. For example, if $E=\bigcup U_i \morph B$ is a cover of a topological space $B$, and $F\morph B$ is a presheaf on $B$, then $E\times_B F$ consists of the values of $F$ on the open cover $\{U_i\}$, $E\times_B E\times_B F$ consists of the values of $F$ on intersections $\{U_i\cap U_j\}$, and so on. The condition for $F$ to be a \emph{sheaf} is evidently equivalent to demanding an equivalence
\[\xymatrix@C=2em{
    F \ar[r]^-\sim
  & \left(E\times_B F \ar@<0.75ex>[r] \ar@<-0.75ex>[r] \right.
  & E\times_B E \times_B F \ar@<1ex>[r] \ar[r] \ar@<-1ex>[r]
  & \left.\dotsb\right)
}\]
of $F$ with its simplicial resolution given by the comonad $E\times_B(-)$. We refer to \cite[\S2]{Hess} for a review of this formalism (and some examples) in the classical setting.

\begin{example}[\protect{\cite{GAGA}}]
\newcommand{\Xan}{X^{\mathrm{an}}}
\newcommand{\XZar}{X^{\mathrm{Zar}}}
\newcommand{\Coh}{\operatorname{Coh}}

Let $X$ be a complex algebraic variety. We can consider $X$ in the analytic topology $\Xan$ with the sheaf $\cH$ of holomorphic functions, or in the Zariski topology\footnote{Here we mean the classical Zariski topology, with no schemy generic points.} $\XZar$ with the sheaf $\cO$ of regular functions. It is not too difficult to show that the forgetful function
\begin{equation}\label{eq:gaga-descent}
(\Xan, \cH) \morph (\XZar, \cO)
\end{equation}
is a faithfully flat map of locally ringed spaces, which is a good notion of ``cover''.

The map \eqref{eq:gaga-descent} induces a functor
\[ \Coh(\cO) \morph \Coh(\cH) \]
between categories of coherent sheaves of modules. The main theorem of \cite{GAGA} is that this is an equivalence of categories when $X$ is projective. We can interpret this as saying that coherent sheaves descend along the cover \eqref{eq:gaga-descent} when $X$ is projective. This is \emph{false} for general $X$, even for $X$ affine.
\end{example}

Descent theory can be formulated using only elementary category theory, and so it is easy to \emph{ask} descent questions. The preceding example shows that \emph{answering} descent questions can involve deep mathematics. It is thus desirable to have very general theorems on when descent holds, which in particular applications may be further simplified to explicit, easily-checkable criteria. The general formalism involves a ``comparison functor'' $k$, and the two basic theorems of general monadic descent theory concern when this functor is fully faithful or an equivalence of categories (we say that \emph{descent} is satisfied in the first case and \emph{effective descent} in the second).

  So far, all this is classical. The rise of derived algebraic geometry and derived stacks has contributed to growing consumer demand for higher descent theory. Following the classical case, this should be formulated in terms of $\infty$-monads. We indicate the challenges in doing so, then explain our approach to surmounting them.

Mathematical theories frequently admit both an \emph{extrinsic} (``choosing coordinates'') as well an \emph{intrinsic} (``coordinate-free'') approach. Typically, the extrinsic approach is useful for carrying out calculations but awkward for developing general theory, while the reverse is true of the intrinsic approach. In abstract homotopy theory/higher category theory\footnote{Opinion is divided on whether or not these terms are synonymous.}, the ``extrinsic'' approach is to ``model'' an $\infty$-category via a ordinary category equipped with additional data specifying the homotopical structure (ideally a simplicial model category). One can then work with the familiar strict morphisms, co/limits, \dots, as long as one makes homotopical corrections along the way (co/fibrant replacements, deriving functors, \dots). In contrast, the ``intrinsic'' approach is to work in a environment where everything is ``fully derived''; as we shall see, an \emph{$\infty$-cosmos} is an extremely robust such environment.

The strategy of working strictly and making homotopical adjustments along the way is extremely effective for a great deal of $\infty$-categorical work (as evidenced by the ubiquity of model categories in the literature). It becomes problematic when working with $\infty$-monads: the equations defining a point-set monad will rarely continue to hold after we make homotopical corrections, thus destroying the strictness which is the point of model categories in the first place. This is compounded when we take iterated composites of a monad. Obviously, this presents a problem for higher descent theory. In particular, while papers such as \cite{Hess} and \cite{AroneChing} have had some success in treating $\infty$-monads and higher descent model-categorically, they must demand fairly stringent hypotheses on the model categories and/or monads involved in order to do so. Although Blumberg-Riehl were able to remove these hypotheses in \cite{BlumbergRiehl}, using the theory of algebraic model categories, control over the $\infty$-category of algebras remained elusive. In view of the preceding discussion, it is natural to move to a fully derived environment in order to treat the foundations of higher descent.

At present, the most comprehensive such environment is that of quasicategories, developed by Joyal and Lurie. Lurie has indeed proven a Barr-Beck theorem in this context \cite[4.7.4.5]{HA}. Subsequently, Riehl and Verity gave a new proof \cite[7.2.7]{RV2}, working in the more general context of \emph{$\infty$-cosmoi}. However, the Barr-Beck theorem only addresses the question of when the monadic comparison functor is an equivalence. As mentioned above, it is also important to know when it is merely fully faithful. The purpose of this paper is to establish this criterion in the $\infty$-categorical setting. We shall deploy the Riehl-Verity framework in order to prove:

{
\renewcommand{\thetheorem}{\ref{thm:main}}
\begin{theorem}
Let $\ds\adjnctn XAfu$ be a homotopy coherent adjunction between $\infty$-categories, inducing a homotopy coherent monad $t=uf$ on $X$ and homotopy coherent comonad $g=fu$ on $A$. Then the comparison functor $A\morph^k \Alg Xt$ to the $\infty$-category of homotopy coherent $t$-algebras is fully faithful if and only if every object of $A$ is $g$-cocomplete, i.e.\ weakly equivalent to the geometric realization of the simplicial resolution given by $g$.
\end{theorem}
}

We now extol the virtues of the Riehl-Verity framework. Classically, it has proven fruitful to develop category theory by working in a nice (behaving like $\Cat$) 2-category. Thus one trades explicit constructions for 2-universal properties. The advantage of this method is that it applies simultaneously to develop the theory of more general categorical structures, such as enriched, internal, or parametrized categories. This approach is often referred to as \emph{formal category theory}, e.g.\ in \cite{Gray}; one could succinctly describe the Riehl-Verity project as \emph{formal $\infty$-category theory}, and an $\infty$-cosmos as an $(\infty,2)$-category behaving like that of $(\infty,1)$-categories (or even $(\infty,n)$-categories).

One thus expects to characterize $\infty$-categorical constructions via $(\infty,2)$- (modelled as simplicially enriched) universal properties. Remarkably, though, the majority of the theory takes place in the \emph{homotopy 2-category}, and so these universal properties are close or identical to those we'd find in the classical case. Sufficiently slick classical proofs\footnote{The task of converting a down-to-earth classical proof into a sophisticated (2-categorical) one---suitable for interpretation in an $\infty$-cosmos---is not necessarily a trivial one.} can thus be transported nearly word-for-word into the $\infty$-categorical context. Indeed, once we get the definitions out of the way, the reader will note we make scarce explicit reference in \S\ref{sec:results} to the definitions of $\infty$-cosmoi.

We now turn to the outline of the paper.

In \S\ref{sec:background} we review the definitional framework and results of Riehl-Verity that we need; this section is expository and discursive, and only sketches of proofs are to be found therein. Readers familiar with their work may skip to \S\ref{sec:results}, which begins with a notational review for the convenience of those who do so. Our results are contained in \S\ref{sec:results}; we explain how to interpret Theorem \ref{thm:main} in an $\infty$-cosmic environment, and prove it. We then indicate some applications to descent problems, including descent spectral sequences.

Finally, we state our position on the most controversial question in the whole of $\infty$-cosmology: how to spell the plural of $\infty$-cosmos. The reader will already have observed that we adhere to the convention of the pioneering $\infty$-cosmologists. We have nothing further to say on the matter, except to affirm that, when we go out for a ramble on a cold day, we do indeed carry supplies of hot coffee with us in \emph{thermoi} \cite{Johnstone}.

\subsection{Acknowledgments}
I would like to thank my advisor, Andrew Blumberg, for suggesting this project and for his guidance and patience throughout. The inspiring question was asked by David Nadler. We are grateful to Emily Riehl for several helpful conversations and clarifications. Some of these conversations took place at the Workshop on Homotopy Type Theory and Univalent Foundations of Mathematics; we thank the organizers for putting the workshop together and for providing travel support. Finally, the typesetting of this paper has benefitted from Dominic Verity's \TeX nical virtuosity.

\section{Background}

\label{sec:background}
Here we review the necessary parts of the prior work of Riehl and Verity. Full details are available in \cite{RV1}, \cite{RV2}, and \cite{RV4}; we recommend \cite{RV0} for a rapid overview. In \S\ref{sub:infty-cosmoi} we introduce the fundamental notions of an $\infty$-cosmos and its homotopy 2-category; this is the setting in which the rest of our work takes place. In \S\ref{sub:homotopy-coherent-adjunctions} we define homotopy coherent/commutative adjunctions and monads, and recall the descriptions of the categories $\Adj$ and $\Mnd$ which corepresent these. Comma $\infty$-categories, which are key to the ``model independence'' of Riehl-Verity's results, are reviewed in \S\ref{subs:commas}. Limits and colimits inside $\infty$-categories are discussed in \S\ref{sub:lims-colims}. In \S\ref{sub:weighted-limits}, we review the enriched-categorical notion of \emph{weighted limits} and discuss their use in the $\infty$-cosmic context, which is simplicially enriched. Finally, \S\ref{sub:monadic-adjunction} shows how to construct the various $\infty$-categories and functors relevant to discussions of monadicity and descent.

\subsection{\texorpdfstring{$\infty$}{Infinity}-cosmoi}
\label{sub:infty-cosmoi}

Informally speaking, an $\infty$-cosmos is a presentation of an $(\infty,2)$-category which is sufficiently well-behaved for us do ``formal $\infty$-category theory'' (à la \cite{Gray}) inside it. (The name is meant to evoke Street, not Bénabou, cosmoi.) The definition is reminiscent of the properties enjoyed by fibrant objects in any model category enriched (c.f.\ \cite[\S A.3.2]{HTT}) over the Joyal model structure on $\sSet$, and indeed these are examples \cite[2.2.1]{RV4}. The reference for this section is \cite[\S2]{RV4}.

\begin{definition}[$\infty$-cosmos]
Let $\cE$ be a simplicially enriched category, equipped with two distinguished classes of 1-cells: the \emph{weak equivalences}, denoted $\acyc$, and the \emph{isofibrations}, denoted $\fibr$. A map which is both a weak equivalence and an isofibration will be called an \emph{acyclic fibration} and denoted $\afibr$. We assume that weak equivalences satisfy the 2-of-6 property, that isofibrations are closed under composition, and that all isomorphisms are acyclic fibrations.

We shall say that $\cE$ is an \emph{$\infty$-cosmos} it it satisfies the following axioms:

\begin{enumerate}
\item (completeness) as a simplicially enriched category, $\cE$ possesses a terminal object $\term$, cotensors $E^J$ of objects $E$ by all simplicial sets $J$, and pullbacks of isofibrations along any functor;

\item (fibrancy) all of the maps $E\fibr\term$ are isofibrations;

\item (pullback stability) isofibrations and acyclic fibrations are stable under pullback along any functor;

\item (SM7) if $E\fibr^p B$ is an isofibration in $\cE$ and $I\monic^i J$ is an inclusion of simplicial sets then the Leibniz cotensor $i \lcot p\colon E^J \to E^I \times_{B^I} B^J$ is an isofibration, and further an acyclic fibration whenever $p$ is an acyclic fibration in $\cE$ or $i$ is an acyclic cofibration in $\sSet_{\text{Joyal}}$; and

\item (cofibrancy) all objects enjoy the left lifting property with respect to all acyclic fibrations in $\cE$.
\end{enumerate}

We will also require an $\infty$-cosmos to have limits of transfinite towers of isofibrations, and for isofibrations to be stable by retracts. We write $\map(A,B)\in\sSet$ for the mapping space between two objects $A,\, B$ of $\cE$.

For psychological reasons, we refer to the objects of $\cE$ as \emph{$\infty$-categories} and its arrows as \emph{functors}.
\end{definition}

\begin{remark}
The axioms can be made stronger or weaker, depending on what one wants to prove. More fundamental is the \emph{style} of arguing about $\infty$-categories: one can imagine working with quasicategories throughout, and the axioms record those properties of quasicategories we use (which turn out to be satisfied much more generally). For example, in \cite{RV4} the axioms only require cotensors by simplicial sets with finitely many non-degenerate simplices; and the ability to take limits of transfinite towers of isofibrations is absent altogether. Our ``infinitary'' assumptions are necessary for the constructions in \S\ref{sub:monadic-adjunction} and \S\ref{sub:main}.
\end{remark}

\begin{remark}
Our assumption that all objects are cofibrant has the crucial consequence that the mapping spaces $\map(A,B)$ are actually \emph{quasicategories}. We refer to \cite[\S2.2]{RV1} for a review of quasicategories. One can get by by merely assuming that every object of $\cE$ has a cofibrant replacement (as in \cite{RV4}, for example); we have chosen not to do so in order to simplify the exposition.
\end{remark}

\begin{definition}
If $A$ is an $\infty$-category in an $\infty$-cosmos $\cE$, the the \emph{underlying quasicategory of $A$} is $\map_\cE(\term, A)$. We define \emph{objects} and \emph{maps} in abstract $\infty$-categories in terms of their underlying quasicategories.
\end{definition}

\begin{example}
In \cite[\S2.2]{RV4}, Riehl and Verity present several ways to easily produce examples of $\infty$-cosmoi. Chief among these examples are:
\begin{itemize}
  \item $\Cat$, the $\infty$-cosmos of ordinary categories. Weak equivalences are equivalences of categories, and isofibrations are functors with the right lifting property with respect to $\{\bullet\}\mono\{\bullet\isom\star\}$.

  \item $\qCat$, the $\infty$-cosmos of quasicategories. Weak equivalences and isofibrations are as usual.

  \item The $\infty$-cosmos of $\theta_n$-spaces, a model of $(\infty,n)$-categories.
\end{itemize}
\end{example}

\begin{example}[\protect{\cite[2.1.11]{RV4}}]
If $A$ is an $\infty$-category in an $\infty$-cosmos $\cE$, we let $\cE/A$ denote the full simplicial subcategory of the usual simplicial slice category spanned by the \emph{isofibrations} $B \fibr A$. This is again an $\infty$-cosmos, called the \emph{slice $\infty$-cosmos} over $A$. Thus the $\infty$-cosmic framework captures parametrized $\infty$-category theory.
\end{example}

With the $\infty$-cosmic framework in hand, Riehl and Verity are able to rederive a great deal of the theory of $\infty$-categories. Their proofs are ``formal'' in nature---in contrast to the combinatorial arguments of \cite{HTT}---and thus permit arguments very close to the classical case. Moreover, as the above examples indicate, their work is not limited to developing the category theory of $(\infty,1)$-categories: it simultaneously applies to develop the category theory of $(\infty,n)$-categories and recapture that of ordinary categories.

However, the import of their work is not merely that $\infty$-cosmoi provide a robust environment in which to develop the category theory of $\infty$-categories. They also show (somewhat unexpectedly) that a much simpler structure suffices for much of this development.

\begin{definition}
The \emph{homotopy 2-category} of an $\infty$-cosmos $\cE$ is the (strict) 2-category $\htc\cE$ with the same underlying category as $\cE$, but with hom-categories $\hom(E, F)$ given by
\[ \hom(E, F) \defeq h(\map(E,F)) \]
for $E,\,F\in\cE$. Here $h$ sends a quasicategory (or simplicial set) to its homotopy category.
\end{definition}

\begin{remark}
When we drop the assumption that all objects in $\cE$ are cofibrant, the 2-category just defined is notated $h_*\cE$, and the correct definition of $\htc\cE$ is the full subcategory of $h_*\cE$ spanned by the (images of) cofibrant objects of $\cE$.
\end{remark}

Recall that a 1-cell $A\morph^f B$ in a 2-category $\cC$ is an \emph{equivalence} if there is a 1-cell $B\morph^g A$ and isomorphic 2-cells $1_A \twomorph^{\sim\,} gf$ and $fg \twomorph^{\sim\,} 1_B$. The following proposition is one of the first indications that $\htc\cE$ remembers enough information about $\cE$ to develop the category theory of its objects.
\begin{proposition}[\protect{\cite[3.1.8]{RV4}}]
A functor $A\morph B$ is a weak equivalence in the $\infty$-cosmos $\cE$ if and only if it is an equivalence in the homotopy 2-category $\htc\cE$.
\end{proposition}
For this reason, we will sometimes write $A=B$ to mean that there exists a weak equivalence $A\acyc B$ in $\cE$ (or, if there is an obvious map $A\morph B$ in play, that this particular map is a weak equivalence).

\subsection{Homotopy coherent adjunctions}
\label{sub:homotopy-coherent-adjunctions}
The reference for this section is \cite[\S3]{RV2}.

\begin{definition}
Let $\cC$ be a 2-category. An \emph{adjunction} in $\cC$ consists of
\begin{itemize}
\item a pair of objects $X$, $A$ of $\cC$;
\item maps $X \morph^f A$ and $A\morph^u X$;
\item two-cells $1_X \twomorph^\eta uf$ and $fu\twomorph^\epsilon 1_A$;
\item satisfying the triangle identities $\epsilon f\cdot f \eta = 1_f$ and $u\epsilon\cdot \eta u = 1_u$.
\end{itemize}

We call $f$ the \emph{left adjoint}, $u$ the \emph{right adjoint}, $\eta$ the \emph{unit}, and $\epsilon$ the \emph{counit} of the adjunction. We indicate an adjunction by writing $f\adjoint u$, $\ds\adjnctn XAfu$ or $\adjnctn XAfu$.
\end{definition}
\begin{definition}
Let $\cC$ be a 2-category. A \emph{monad} in $\cC$ consists of an object $X\in\cC$ and a monoid $t$ in the monoidal category $\hom(X,X)$.
\end{definition}

When $\cC=\Cat$, these specialize to the usual notions. Since these notions are equationally defined, they are corepresentable, i.e.\ there is a 2-category $\Adj$ (resp.\ $\Mnd$) such that adjunctions (resp.\ monads) in $\cC$ are the same thing as 2-functors $\Adj\to\cC$ (resp.\ $\Mnd\to\cC$). The explicit description of $\Adj$ is due to Schanuel and Street \cite{SchanuelStreet}, of $\Mnd$ to Lawvere \cite{Lawvere}. Before giving the definition, we set some notation.

\begin{definition}
As usual, $\Delta_+$ and $\Delta$ will denote the category of finite linearly ordered sets and the full subcategory of non-empty sets. We shall use the notation $\Delta_\infty$ (respectively $\Delta_{-\infty}$) to denote the subcategory of $\Delta$ consisting of those maps which preserve top (respectively bottom) elements.
\end{definition}

\begin{definition}
The free adjunction is the small 2-category $\Adj$ with two objects $+$ and $-$, with hom-categories given by
\begin{align*}
  \Adj(+,+) &= \Delta_+ & \Adj(-,-) &= \Delta_+^\op\\
  \Adj(-,+) &= \Delta_\infty \isom \Delta_{-\infty}^\op & \Adj(+,-) &= \Delta_{-\infty} \isom \Delta_\infty^\op
\end{align*}
as summarized in the following picture:
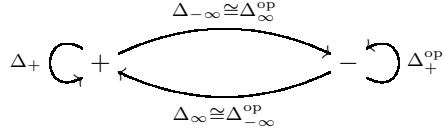
\[
  \xymatrix@C=8em{
    {+}\ar@/^3ex/[r]^{\Delta_{-\infty}\cong\Delta_\infty^\op}
    \ar@(ul,dl)[]_{\Delta_+} &
    {-}\ar@/^3ex/[l]^{\Delta_\infty\cong\Delta_{-\infty}^\op}
    \ar@(dr,ur)[]_{\Delta_+^\op} \\
  } 
\]
We write $+\morph^\Adjf-$ (resp.\ $-\morph^\Adju+$) for the map corresponding to $[0]\in\Delta_{-\infty}$ (resp.\ to $[0]\in\Delta_\infty$).
\end{definition}

\begin{definition}
The free monad is the small 2-category $\Mnd$ which is the full subcategory of $\Adj$ on the object $+$. We write $\Adjt$ for the endomorphism corresponding to $[0]\in\Delta_+$.
\end{definition}

\begin{definition}
Any 2-category gives rise to a simplicially enriched (in fact, quasicategorically enriched) category by identifying the hom-categories with their nerves (this uses the fact that the nerve preserves products). This process is the right adjoint $N_*$ in a 2-adjunction
\[
  \adjnctn{\text{(simplicial categories)}}{\text{(2-categories)}}{h_*}{N_*}
  \vspace{-2ex}
\]
arising from the Quillen adjunction $\ds\adjnctn{\sSet_{\text{Joyal}}}\Cat hN$; we have already made use of $h_*$ in defining the homotopy 2-category $\htc\cE$ of an $\infty$-cosmos $\cE$.

Applying this to the 2-categories $\Adj$ and $\Mnd$, we obtain simplicially enriched categories which we continue to notate $\Adj$ and $\Mnd$. Since $N$ is fully faithful, this conflation is anodyne.
\end{definition}

\begin{remark}
The calculus of string diagrams for 2-categories extends naturally to describe the $n$-arrows of simplicial categories which arise in this way. Riehl and Verity show in \cite{RV2} that when specialized to $\Adj$, this graphical calculus admits a variation---the calculus of ``strictly undulating squiggles''---enabling a simple combinatorial description of the $n$-arrows of $\Adj$ which behaves well with respect to both vertical and horizontal composition. Strikingly, they use this to show that $\Adj$ is cofibrant in the Bergner model structure on simplicial categories \cite{Bergner}, and to work with explicit cellular presentations of $\Adj$.
\end{remark}

\begin{notation}
The symbol $-$ is often used as a placeholder symbol in category theory. To avoid confusion with the object $-$ of $\Adj$, we will use $\slot$ instead. Thus $\Adj(\slot,+)$ is a functor $\Adj^\op\to\sSet$, but $\Adj(-,+)$ is an object of $\sSet$.
\end{notation}

\begin{definition}
Let $\cE$ be an $\infty$-cosmos with homotopy 2-category $\htc\cE$.
\begin{itemize}
  \item A \emph{homotopy coherent adjunction}, or \emph{$\infty$-adjunction}, in $\cE$ is a simplicial functor $\Adj\to\cE$.

  \item A \emph{homotopy commutative adjunction}, or \emph{1-adjunction}, in $\cE$ is a 2-functor $\Adj\to\htc\cE$.

  \item A \emph{homotopy coherent monad}, or \emph{$\infty$-monad}, in $\cE$ is a simplicial functor $\Mnd\to\cE$.

  \item A \emph{homotopy commutative monad}, or \emph{1-monad}, in $\cE$ is a 2-functor $\Adj\to\htc\cE$.
\end{itemize}
\end{definition}

\begin{warning}
When an $\infty$-category $X$ is presented by a 1-category (e.g.\ a simplicial model category) $\cX$, 1-monads on $X$ as defined above must not be confused with ``point-set'' monads on $\cX$. The former are monads on $\Ho(\cX)$; the latter sometimes induce $\infty$-monads on $X$, but we shall make no pre-derived use of them. It does not appear to be possible to give a simple definition of $\infty$-monads on $X$ purely in terms of $\cX$ unless $\cX$ is very special, e.g.\ a simplicial model category in which everything is bifibrant.
\end{warning}

\begin{litcomp}\label{lc:adj}
Let $\ds\smash[t]{\adjnctn \cX\cA FU}$ be a simplicial Quillen adjunction between simplicial model categories $\cX$ and $\cA$. $\cX$ and $\cA$ model $(\infty,1)$-categories $X$ and $A$. For example, to obtain quasicategorical models, we would take homotopy coherent nerves of the subcategories of bifibrant objects: $X \defeq N(\cX_\cf)$ and $A\defeq N(\cA_\cf)$. By \cite[6.2.1]{RV1}, there is an induced $\infty$-adjunction $\ds\adjnctn XAfu$. These functors are obtained by correcting $F$ and $U$ to land in bifibrant objects. For example, if every object of $\cX$ is cofibrant and every object of $\cA$ is fibrant, then no correction is needed.

Work of Dugger, Rezk, Schwede, and Shipley shows that a Quillen adjunction between left proper combinatorial model categories is equivalent to a simplicial Quillen adjunction as above; see \cite[\S A]{BlumbergRiehl} for discussion of this. Thus we again get an induced $\infty$-adjunction between $\infty$-categories.
\end{litcomp}

\begin{litcomp}
Let $X$ be an $\infty$-category modelled by a simplicial model category $\cX$, and let $\bbT$ be a simplicial monad on $\cX$. Under reasonable conditions, the category $\Alg\cX\bbT$ of $\bbT$-algebras is a simplicial model category in such a way that the monadic adjunction $\adjnctn \cX{\Alg\cX\bbT}{F^\bbT}{U^\bbT}$ is simplicial Quillen \cite[\S C]{Hess}. We thus obtain an $\infty$-adjunction out of $X$, and hence an $\infty$-monad $t$ on $X$.
\end{litcomp}

Theorems 4.3.9, 4.3.11, 4.4.11, and 4.4.18 of \cite{RV2} show that every 1-adjunction $\Adj\to\htc\cE$ lifts to an $\infty$-adjunction $\Adj\to\cE$, and moreover that such lifts are unique in a suitable homotopical sense. The proof proceeds by explicit analysis of the combinatorics of such lifting problems, made possible by the squiggle calculus mentioned above. In contrast, it is \emph{not} possible in general to lift 1-monads to $\infty$-monads.

\subsection{Comma \texorpdfstring{$\infty$}{infinity}-categories}
\label{subs:commas}

\begin{definition}
Let $B\morph^f A \comorph^g C$ be a pair of functors in the $\infty$-cosmos $\cE$. The \emph{comma $\infty$-category} $\comma fg$ is defined to be the pullback in the diagram below:
\[\xymatrix{
  \comma fg \pullback \ar@{->>}[d] \ar[r] & A^{\Delta^1} \ar@{->>}[d]^-{\<\cod,\,\dom\>}\\
  C\times B \ar[r]_-{g\times f} & A\times A
}\]
In the case of an identity functor, we write $\comma fA$ instead of $\comma f{\id_A}$.
\end{definition}

\begin{remark}
$\comma fg$ should be thought of as having objects triples $\<b\in B, c\in C, fb\morph^\phi gc\in A\>$. The reason for writing $C\times B$ and $\<\cod,\,\dom\>$ instead of the seemingly more natural $B\times C$ and $\<\dom,\,\cod\>$ in the above diagram is that $\comma fg$ is a ``(left $C$, right $B$)-bimodule''. By this we mean that $\comma fg$ carries a covariant action by $C$ and a contravariant action by $B$, and these commute. This perspective is extremely useful, and is the subject of \cite{RV5}. Some of the proofs in \S\ref{sub:descent} use the 1-categorical (but simplicially enriched) version of this ``calculus of modules'', for which a good reference is \cite[\S\S4.1 and 4.3]{RiehlCHT}.
\end{remark}

\begin{remark}
Comma $\infty$-categories are important for (at least) two reasons. First, general 2-category theory would have us define many notions representably, carrying around ``generalized objects'' $Z\to A$ (since we can't ``look inside'' our $\infty$-categories). Comma categories allow us to dispense with this artifice and work more directly with the $\infty$-category $A$, thus keeping our intuition close to the classical case. Secondly, as we shall see, all the basic notions of category theory can be expressed in terms of commas. Once it is shown that functors of $\infty$-cosmoi preserve commas \cite[2.3.10]{RV5}, it follows that $\infty$-category theory developed in the $\infty$-cosmic framework is ``model independent''. See \cite[\S\S3.6 and 4.5]{RV0} for further discussion.
\end{remark}

\begin{example}[\protect{\cite[4.4.2 and 4.4.3]{RV1}}]
Let $X\morph^f A$ and $A\morph^u X$ be a pair of functors between $\infty$-categories. Then $f\adjoint u$ if and only if $\comma fA = \comma Xu$ in the slice $\infty$-cosmos $\cE/(A\times X)$:
\[\xymatrix@C=1ex{
  \comma fA \ar@{->>}[dr] \ar[rr]^-\sim && \comma Xu \ar@{->>}[dl]
  \\ & A \times X
}\]
\end{example}

\begin{example}
If $a$ and $b$ are objects of an $\infty$-category $A$, the comma $\infty$-category $\comma ab$ is a model of the mapping space between $a$ and $b$ inside $\cE$.
\end{example}

\begin{notation}
If $a$ and $b$ are objects of an $\infty$-category $A$, we write $A(a,b) \defeq \map(\term,\comma ab)$ for the underlying quasicategory of $\comma ab$; by \cite[\S3.2]{RV0}, this is a Kan complex. If a map $a\morph^f b$ is given, we write $A(a,b)_f$ to denote $A(a,b)$ equipped with the basepoint $f$.
\end{notation}

\begin{lemma}\label{lm:cot-comma}
If $J\in\sSet$, then ${\comma fg}^J = \comma {f^J}{g^J}$.
\end{lemma}

\subsection{Limits and colimits in \texorpdfstring{$\infty$}{infinity}-categories}
\label{sub:lims-colims}

The reference for this section is \cite[\S5]{RV1}.

\begin{definition}
An \emph{absolute left lifting diagram} in a 2-category consists of the data
\begin{equation}
  \label{eq:alld}
  \vcenter{\xymatrix{
    \ar@{}[dr]|(0.7){\Uparrow\lambda} & C \ar[d]^-\psi\\
    A \ar[r]_-\phi \ar[ur]^-\ell & B
  }}
\end{equation}
inducing unique factorizations of 2-cells:
\[
  \vcenter{\xymatrix{
    Z \ar[r]^-q \ar[d]_-p \ar@{}[dr]|(0.5){\Uparrow\zeta} & C \ar[d]^-\psi\\
    A \ar[r]_-\phi & B
  }}
  \quad=\quad
  \vcenter{\xymatrix{
    Z \ar[r]^-q \ar[d]_-p
    \ar@{}[dr]|(0.3){\exists!\Uparrow} \ar@{}[dr]|(0.7){\Uparrow\lambda}
    &
    C \ar[d]^-\psi\\
    A \ar[r]_-\phi \ar[ur]|{\hole\ell\hole} & B.
  }}
\]
\end{definition}
The pasting operation can be broken down as
\[\xymatrix@C=0em{
  \hom(Z, C)(\ell p, q) \ar[rr]^-{\text{paste with $\lambda$}} \ar[dr]_-{\text{whisker with }\psi\qquad} && \hom(Z, B)(\phi p, \psi q)\\
  & \hom(Z, B)(\psi \ell p, \psi q) \ar[ur]_-{\quad\qquad\text{precompose vertically with }\lambda}
}\]
and the definition is demanding that ``paste with $\lambda$'' be a bijection for all spans $A\comorph^p Z\morph^q C$.

It will be useful to characterize absolute lifting diagrams in terms of comma categories rather than a test object $Z$. In fact, \eqref{eq:alld} is an absolute left lifting diagram if and only if the map $\comma\ell F \to \comma \phi\psi$ induced by $\lambda$ is an equivalence \cite[5.1.3]{RV1}.

Finally, we note that \eqref{eq:alld} is an absolute left lifting diagram if and only if it has that property pointwise, that is, for all $a\in A$ \cite[6.1.8]{RV1}.

\begin{definition}
Let $J\in\sSet$. We say that an $\infty$-category $E$ \emph{admits colimits of a family of diagrams $D\morph^d E^J$ of shape $J$} if there is an absolute left lifting diagram
\[\xymatrix{
  \ar@{}[dr]|(0.7){\Uparrow\lambda} & E \ar[d]^-c\\
  D \ar[ur]^-{\colim} \ar[r]_-d & E^J
}\]
in $\htc\cE$; here $c$ is the constant map. In this case we call $\lambda$ a \emph{colimiting cone}.
\end{definition}

The definition asks for the existence of a functor $D\morph E$ and a 2-cell $\lambda$ satisfying certain universal properties. In our work in \S\ref{sub:main}, the functor and 2-cell will always exist: the question will be \emph{whether} they define an absolute left lifting diagram. As mentioned above, this is the case if and only if the map $\comma\colim E \to \comma dc$ induced by $\lambda$ is an equivalence.

\begin{remark}
When the $\infty$-cosmos $\cE$ is cartesian closed, one can give a completely analogous definition of colimits for shapes $J\in\cE$. If $\cE$ is not cartesian closed, a different approach must be used; see \cite{RV5}. We shall only require diagram shapes given by simplicial sets.
\end{remark}

\subsection{Weighted limits}
\label{sub:weighted-limits}

The last section discussed limits \emph{in} $\infty$-categories; we will also require limits \emph{of} $\infty$-categories. In \S\ref{sub:monadic-adjunction} this will be employed to tame the zoo of $\infty$-categories unleashed by an adjunction, by characterizing them by universal properties. The first half of this section is our telling of a standard story; the reference for the second half is \cite[\S5.2]{RV2}.

Being in the context of simplicially enriched categories imposes enriched category theory on us. Limits are one area where very different behavior arises in the enriched world than for ordinary categories: ordinary limits still make sense in the enriched case, but are woefully inadequate. The enriched context demands we consider \emph{weighted limits}, a notion we suggest some intuition for before giving the precise definition.

Let $\cV$ be a \emph{Bénabou} cosmos: a bicomplete closed symmetric monoidal category (``a category suitable for enriching over''). We shall only need $\cV=\sSet$, but the theory is perfectly general. Let $\cA\morph^T\cC$ be a $\cV$-functor between $\cV$-categories $\cA$ and $\cC$, with $\cA$ small. Remember that the \emph{limit} $\projlim T$ of $T$ is defined by requiring it to represent cones over $T$; that is, we have a natural correspondence between
\[
  Z \morph^{\phi} \projlim T
  \qquad\text{and}\qquad
  \vcenter{\xymatrix{
    & Z \ar[dl]_-{\phi_a} \ar[dr]^-{\phi_{a'}} \\
    Ta \ar@<1ex>[r] & \ar@<1ex>[l] \cdots \ar@<-1ex>[r] & Ta'. \ar@<-1ex>[l]
  }}
\]
In the enriched context, we can demand richer structure $Wa\morph^{\phi_a}\cC(Z, Ta)\in\cV$ than just specifying a single map $\phi_a$, and we define $\wLim WT\cA$, the \emph{limit of $T$ weighted by $W$}, by demanding a natural correspondence between
\[
  Z \morph^{\phi} \wLim WT\cA
  \qquad\text{and}\qquad
  \vcenter{\xymatrix{
    & Z \ar[dl]|-{\hole\phi_a(Wa)\hole} \ar[dr]|-{\hole\phi_{a'}(Wa')\hole} \\
    Ta \ar@<1ex>[r] & \ar@<1ex>[l] \cdots \ar@<-1ex>[r] & Ta'. \ar@<-1ex>[l]
  }}
\]
More precisely, let $\cA\morph^W \cV$ be a $\cV$-functor, which we call the \emph{weight}. The weighted limit $\wLim WT\cA$ is defined by the universal property
\begin{equation}\label{eq:wlims}
  \cC(Z, \wLim WT\cA) = \cV^\cA(W, \cC(Z, T(\slot))).
\end{equation}
Important cases include $\wLim {\text{(constant at object }V\in\cV)}T\cA = \left(\projlim T\right)^V = \projlim T(\slot)^V$, a cotensor of the ordinary limit, and $\wLim {\cA(a,-)}T\cA = T(a)$ (the latter, as usual, is more or less the Yoneda lemma). Importantly, $\wLim\slot T\cA\colon (\cV^\cA)^\op \to \cC$ is a right adjoint, and so takes colimits of weights to limits of weighted limits. Combined with the two cases just mentioned, this gives the end formula
\[ \wLim WT\cA = \int_{a\in\cA} Ta^{Wa} \]
which in particular shows that having all weighted limits is equivalent to having all ordinary limits and all cotensors over $\cV$.

In general, an $\infty$-cosmos will not have all weighted limits. However, there is a conceptually elegant description of the weighted limits which do exist. If $\cA$ is a small simplicial category, say that a natural transformation in $\sSet^\cA$ is a \emph{projective cofibration} if it has the left lifting property with respect to level acyclic fibrations. The projective cofibrations are evidently the closure of the set
\[
  \{ \bdry\Delta^n \times \cA(a, \slot) \monic \Delta^n \times \cA(a, \slot) \st n\ge 0,\, a\in\cA \}
\]
of \emph{projective cells} in the Galois correspondence defined by left/right lifting properties. In particular, a natural transformation is a projective cofibration if and only if it is a retract of a transfinite composite of pushouts of projective cells.

\begin{proposition}[\protect{\cite[5.2.4]{RV2}}]\label{prop:proj-cof-lims}
An $\infty$-cosmos has all limits weighted by projective cofibrant weights.
\end{proposition}

Indeed, $\wLim \slot T\cE$ turns all the types of colimits used to build projective cofibrations from projective cells into types of limits which are guaranteed to exist by the $\infty$-cosmos axioms (we added some axioms for precisely this purpose). We are thus reduced to showing that limits weighted by projective cells exist and are isofibrations; but this follows immediately from the completeness and SM7 axioms for an $\infty$-cosmos.

\begin{warning}
If $L$ is an object of $\cE$ is defined by a weighted limit, and thus satisfying a $\sSet$-enriched universal property, it is generally not true that the image of $L$ in $\htc\cE$ will have the analogous $\Cat$-enriched universal property. However, it may satisfy a weaker uniqueness condition guaranteeing its uniqueness up to isomorphism, although not up to automorphisms. Let $L\in\cE$, and let $W\twomorph^\phi \cC(L,T(\slot))$ be a weighted cone. We say that $\phi$ displays $L$ as a \emph{weak 2-limit of $T$ weighted by $W$} if the induced functors
\[\cC(Z,L) \morph \cV^\cA(W, \cC(Z,T(\slot))),\]
rather than being equalities as in \eqref{eq:wlims}, are \emph{smothering}: surjective on objects, full, and conservative. For example, comma categories are weak 2-limits in this sense \cite[3.3.18]{RV1}, a fact which we will use in the proof of Lemma \ref{lm:prt-adj}.

Since the above properties can be given by right lifting properties, it follows that fibres of a smothering functor, while not necessarily \emph{contractible}, are at least (non-empty) connected groupoids (i.e.\ classifying spaces of discrete groups). For further details, see \cite[\S3.3]{RV1}.
\end{warning}

\subsection{Monadic adjunction}
\label{sub:monadic-adjunction}

Suppose given a homotopy coherent monad $t$ on an $\infty$-category $X$. In this section we construct the $\infty$-category $\Alg Xt$ of homotopy coherent $t$-algebras, as well as the monadic adjunction $\ds\smash[b]{\adjnctn X{\Alg Xt}{f^t}{u^t}}$. When $t$ arises from a homotopy coherent adjunction $\ds\adjnctn XAfu$, we construct the comparison functor $A\morph^k \Alg Xt$ which will be the subject of \S\ref{sec:results}. The reference for this section is \cite[\S\S6 and 7]{RV2}.

Denote the corepresentable functors by $\Adj^\pm = \Adj(\pm,\slot)$. Write $\ds\adjnctn {\sSet^\Mnd}{\sSet^\Adj}\lan\res$ for the left Kan extension $\adjoint$ restriction adjunction arising from the inclusion of $\Mnd$ into $\Adj$.

\begin{definition}
Let $t$ be a homotopy coherent monad on $X$, given by a simplicial functor $\Mnd\morph^H\cE$. The \emph{$\infty$-category $\Alg Xt$ of homotopy coherent $t$-algebras} (or \emph{$\infty$-$t$-algebras}) is defined by the weighted limit
\[ \Alg Xt \defeq \wLim{\res\Adj^-}H\Mnd. \]
\end{definition}

This is legitimate by Proposition \ref{prop:proj-cof-lims} and \cite[6.1.8]{RV2}. Since $\Mnd\mono\Adj$ is fully faithful, we have $\wLim{\res\Adj^+}H\Mnd=H(+)=X$. The monadic adjunction $\adjnctn X{\Alg Xt}{f^t}{u^t}$ is defined to be coclassified by $\wLim{\res\Adj^\slot}H\Mnd\colon\Adj\to\cE$.

\begin{litcomp}
Our homotopy coherent $t$-algebras correspond to the strictly $\bbT$-complete objects of \cite[4.14]{Hess}. See also \cite[4.20]{Hess}.
\end{litcomp}

Now suppose that $t$ comes from a homotopy coherent adjunction $\ds\smash{\adjnctn XAfu}\vspace{1ex}$ in $\cE$, coclassified by a simplicial functor $\Adj\morph^T\cE$. An inspection of universal properties shows that
\[ \wLim{\lan\res\Adj^-}T\Adj = \wLim{\res\Adj^-}{\res T}\Mnd = \Alg Xt, \]
so we may take all weighted limits over $\Adj$.

\begin{definition}
The \emph{comparison functor} $A\morph^k \Alg Xt$ is defined by requiring the diagram of $\infty$-categories on the right to be induced by the diagram of weights on the left.
\[
  \vcenter{\xymatrix{
    \lan\res\Adj^- \ar@<1ex>[r] \ar@{}[r]|-\bot \ar@{ >->}[d]_-{\text{counit}} & \Adj^+ \ar@{=}[d] \ar@<1ex>[l]\\
    \Adj^- \ar@<1ex>[r]^-{\Adj^\Adjf} \ar@{}[r]|-\bot & \Adj^+ \ar@<1ex>[l]^-{\Adj^\Adju}
  }}
  \quad\rightsquigarrow\quad
  \vcenter{\xymatrix{
    X \ar@{=}[d] \ar@{}[r]|{\bot} \ar@<1ex>[r]^f & A \ar@<1ex>[l]^u \ar[d]^-k\\
    X \ar@{}[r]|{\bot} \ar@<1ex>[r]^{f^t} & \Alg Xt \ar@<1ex>[l]^{u^t}
  }}
\]
\end{definition}

That is, $k$ is induced by the counit of the $\lan\adjoint\res$ adjunction, valued at $\Adj^-$. \cite[7.1.5]{RV2} shows that $\lan\res\Adj^-$ is the subfunctor of $\Adj^-$ consisting of maps which factor through $+$ (and the counit is the inclusion); in particular, $\lan\res\Adj^-(+) = \Delta_\infty$ and $\lan\res\Adj^-(-)=\spc$.

\section{Comparison and cocompletion}
\label{sec:results}

In this section we prove the main theorem. Background on fully faithful functors appears in \S\ref{sub:ff}. In \S\ref{sub:main} we state and prove our main result, characterizing when the comparison functor induced by a monad is fully faithful in terms of a ``cocomplete'' criterion. Applications to descent, including descent spectral sequences, are discussed in \S\ref{sub:descent}.

We begin by establishing the notation to be used throughout this section, and reviewing that which was introduced in \S\ref{sec:background}. Recall that the ``walking adjunction'' is denoted by $\Adj$, the ``walking monad'' by $\Mnd$; we write $\ds\adjnctn {\sSet^\Mnd}{\sSet^\Adj}\lan\res$ for the resulting left Kan extension $\adjoint$ restriction adjunction. Corepresentable functors are written $\Adj^\pm\defeq\Adj(\pm,\slot)$.

Fix once and for all an $\infty$-cosmos $\cE$ and homotopy coherent adjunction $\Adj\morph^T\cE$. We write $\ds\adjnctn XAfu$ for the image of $\ds\adjnctn +-\Adjf\Adju$ in $\cE$. Let $t=uf$ (resp.\ $g=fu$) be the (co)monad induced on $X$ (resp.\ $A$). The monadic adjunction is denoted $\ds\smash[t]{\adjnctn X{\Alg Xt}{f^t}{u^t}}$, the comparison functor $A\morph^k \Alg Xt$, and the descent comonad $g^t=f^tu^t$. All this is summarized in the following picture.
\[\xymatrix{
  X \ar@(u,l)[]_t \ar@{=}[d] \ar@<1ex>[r]^-f \ar@{}[r]|-\bot
  &
  A \ar@(u,r)[]^g \ar@<1ex>[l]^-u \ar[d]^-k
  \\
  X \ar@<1ex>[r]^-{f^t} \ar@{}[r]|-\bot
  &
  \Alg Xt \ar@(d,r)_-{g^t} \ar@<1ex>[l]^-{u^t}
}\]
Recall that this diagram is obtained by applying $\wLim\slot T\Adj$ to the following diagram of weights:
\[\xymatrix{
  \lan\res\Adj^- \ar@<1ex>[r] \ar@{}[r]|-\bot \ar@{ >->}[d]_-{\text{counit}} & \Adj^+ \ar@{=}[d] \ar@<1ex>[l]\\
  \Adj^- \ar@<1ex>[r]^-{\Adj^\Adjf} \ar@{}[r]|-\bot & \Adj^+ \ar@<1ex>[l]^-{\Adj^\Adju}
}\]

\begin{litcomp}
The algebraic model category approach of \cite{BlumbergRiehl} provides perhaps the closest link between model-categorical input and $\infty$-categorical output. Let $\ds\adjnctn \cX\cA FU$ be a simplicial Quillen adjunction between cofibrantly generated\footnote{Consult Aside 3.4 and the remarks after Theorem 3.3 of \cite{BlumbergRiehl} for precisely what we mean by this term here.} simplicial model categories, inducing an $\infty$-adjunction $\ds\adjnctn XAfu$ between the $\infty$-categories $X=N(\cX_\cf)$ and $A=N(\cA_\cf)$. By \cite[6.1]{BlumbergRiehl}, there is a simplicially enriched fibrant replacement monad $\bbR=(R,r,\mu)$ on $\cA$ and a simplicially enriched cofibrant replacement comonad $\bbQ=(Q,q,\nu)$ on $\cX$; thus $f=N(RF|_{\cX_\cf})$ and $u=N(QU|_{\cA_\cf})$. Let $T=QURF$ and $G=RFQU$, which model the $\infty$-monad $t=uf=N(T|_{\cX_\cf})$ on $X$ and the $\infty$-comonad $g=fu=N(G|_{\cA_\cf})$ on $A$. By \cite[6.3]{BlumbergRiehl}, there are \emph{point-set level} simplicially enriched resolutions
\vspace{0.25em}
\[
  \vspace{0.5em}
  \vcenter{\xymatrix@=1.75em{
    Q \ar[r]|-\zeta & TQ \ar@<1ex>[r]|-\zeta \ar@<-1ex>[r] & T^2Q \ar[l] \ar@<2ex>[r]|-\zeta \ar[r] \ar@<-2ex>[r] & \ar@<1ex>[l] \ar@<-1ex>[l] \cdots\\
    R & GR \ar[l]|-\xi \ar[r] & G^2R \ar@<-1ex>[l]|-\xi \ar@<1ex>[l] \ar@<1ex>[r] \ar@<-1ex>[r] & \ar@<2ex>[l] \ar[l] \ar@<-2ex>[l]|-\xi \cdots\\
    RFQ \ar@<-0.5ex>[r] & \ar@<-0.5ex>[l] GRFQ \ar@<0.5ex>[r] \ar@<-1.5ex>[r] & \ar@<-1.5ex>[l] \ar@<0.5ex>[l] G^2RFQ \ar@<1.5ex>[r] \ar@<-0.5ex>[r] \ar@<-2.5ex>[r] & \ar@<1.5ex>[l] \ar@<-0.5ex>[l] \ar@<-2.5ex>[l] \cdots\\
    QUR \ar@<0.5ex>[r] & \ar@<0.5ex>[l] TQUR \ar@<1.5ex>[r] \ar@<-0.5ex>[r] & \ar@<-0.5ex>[l] \ar@<1.5ex>[l] T^2QUR \ar@<2.5ex>[r] \ar@<0.5ex>[r] \ar@<-1.5ex>[r] & \ar@<2.5ex>[l] \ar@<0.5ex>[l] \ar@<-1.5ex>[l] \cdots
  }}
  \quad\text{presenting}\quad
  \vcenter{\xymatrix@=1.75em{
    1 \ar[r]|-\eta & t \ar@<1ex>[r]|-\eta \ar@<-1ex>[r] & t^2 \ar[l] \ar@<2ex>[r]|-\eta \ar[r] \ar@<-2ex>[r] & \ar@<1ex>[l] \ar@<-1ex>[l] \cdots\\
    1 & \ar[l]|-\epsilon g \ar[r] & \ar@<-1ex>[l]|-\epsilon \ar@<1ex>[l] g^2 \ar@<1ex>[r] \ar@<-1ex>[r] & \ar@<2ex>[l] \ar[l] \ar@<-2ex>[l]|-\epsilon \cdots\\
    f \ar@<-0.5ex>[r] & gf \ar@<-0.5ex>[l] \ar@<0.5ex>[r] \ar@<-1.5ex>[r] & \ar@<-1.5ex>[l] \ar@<0.5ex>[l] g^2f \ar@<1.5ex>[r] \ar@<-0.5ex>[r] \ar@<-2.5ex>[r] & \ar@<1.5ex>[l] \ar@<-0.5ex>[l] \ar@<-2.5ex>[l] \cdots\\
    u \ar@<0.5ex>[r] & tu \ar@<1.5ex>[r] \ar@<-0.5ex>[r] \ar@<0.5ex>[l] & \ar@<-0.5ex>[l] \ar@<1.5ex>[l] t^2u \ar@<2.5ex>[r] \ar@<0.5ex>[r] \ar@<-1.5ex>[r] & \ar@<2.5ex>[l] \ar@<0.5ex>[l] \ar@<-1.5ex>[l] \cdots
  }}
\]
at the level of $\infty$-categories. Here the unit $\eta$ and counit $\epsilon$ of the $\infty$-adjunction $f\adjoint u$ are modelled on the point-set level, in terms of the unit $\mathring\eta$ and counit $\mathring\epsilon$ of the point-set adjunction $F\adjoint U$, by
\[\xymatrix{
  \zeta\from Q \ar[r]^-\nu & Q^2 \ar[r]^-{Q\mathring\eta} & QUFQ \ar[r]^-{QUr} & QURFQ
}\]
and
\[\xymatrix{
  \xi\from RFQUR \ar[r]^-{RFq} & RFUR \ar[r]^-{R\mathring\epsilon} & R^2 \ar[r]^-\mu & R.
}\]
We will summon these assumptions and notations with the phrase, ``suppose given model-categorical input''.
\end{litcomp}

\subsection{Fully faithful functors}\label{sub:ff}

We recall the $\infty$-cosmic definition of ``fully faithful'', and demonstrate some elementary facts about it.

\begin{definition}\label{def:ff}
Let $A\morph^k B$ be a functor in $\cE$, and let $a\in A$. We say that $k$ is \emph{fully faithful on maps out of $a$} if the induced functor $\comma aA\to\comma{ka}k$ is an equivalence. We say that $k$ is \emph{fully faithful} if this holds for all $a\in A$; by \cite[6.1.8]{RV1}, this is equivalent to demanding that $\comma AA\to\comma kk$ be an equivalence.
\end{definition}

These are equivalent to asking that
\[
  \vcenter{\xymatrix{
    & A \ar[d]^-k\\
    \term \ar[ur]^-a \ar[r]_-{ka} & B
  }}
  \qquad\text{or}\qquad
  \vcenter{\xymatrix{
    & A \ar[d]^-k\\
    A \ar@{=}[ur] \ar[r]_-k & B
  }}
\]
be absolute left lifting diagrams. Explicitly, the latter means that $k$ induces a bijection between 2-cells $\twoCell ZApq{}$ and $\twoCell ZB{kp}{kq}{}$ for every parallel pair $\smash{\dimorph ZApq}$ in $\htc\cE$ (i.e., $k$ is ``representably fully faithful'').

\begin{litcomp}
In the $\infty$-cosmos of quasicategories, we recapture Lurie's definition \cite[1.2.10.1]{HTT}. In the $\infty$-cosmos of simplicial categories, we recapture the DK-embeddings (the simplicial functors which are locally weak equivalences of simplicial sets).
\end{litcomp}

The results in the remainder of this subsection are not needed in the sequel, but may be of independent interest. (They appeared in earlier attempts to prove Theorem \ref{thm:main}.)

\begin{lemma}\label{lm:cot-ff}
If $k$ is fully faithful, then so is $k^J$ for any $J\in\sSet$.
\end{lemma}
\begin{proof}
Apply Lemma \ref{lm:cot-comma} and \cite[1.2.7.3]{HTT}.
\end{proof}

\begin{lemma}\label{thm:ff-refl}
A fully faithful functor reflects colimits.
\end{lemma}
\begin{proof}
Consider a diagram of the form
\begin{equation}
  \label{eq:ff-refl}
  \vcenter{\xymatrix{
    \ar@{}[dr]|(.7){\Uparrow\lambda} & A \ar[d]^-c \ar[r]^-k & B \ar[d]^-c\\
      K \ar[ur]^-L \ar[r]_-D & A^J \ar[r]_-{k^J} & B^J
  }}
\end{equation}
in which $k$ is fully faithful and $J\in\sSet$.  The statement means that the triangle is an absolute left lifting diagram whenever the composite diagram is. Consider the commutative diagram
\[\xymatrix{
  \comma LA \ar[r]^-\lambda \ar[d]_-k & \comma Dc \ar[d]^-{k^J}\\
  \comma{kL}k \ar[r]_-{k^J \lambda} & \comma{k^J D}{ck}.
}\]
The vertical arrows are equivalences by assumption and Lemma \ref{lm:cot-ff}. If the composite diagram in \eqref{eq:ff-refl} is an absolute left lifting diagram, then the bottom horizontal arrow is an equivalence. In this case the top horizontal arrow must be an equivalence, which is precisely the condition for the triangle of \eqref{eq:ff-refl} to be an absolute left lifting diagram.
\end{proof}

\subsection{Proof of main result}
\label{sub:main}

In order to state our main result, we require some further notation. We define the maps in the diagram to the right by requiring it to be induced by the diagram of weights on the left, in which the horizontal arrows are given by composition and the vertical arrows by inclusions.
\[
  \vcenter{\xymatrix{
    \lan\res\Adj^-\times\spc \ar[r] \ar[d] & \lan\res\Adj^- \ar[d]\\
    \Adj^-\times\spc \ar[ur] \ar[r] & \Adj^-
  }}
  \rightsquigarrow\quad
  \vcenter{\xymatrix{
    A \ar[r]^-{\resol} \ar[d]_-k & A^\spc \ar[d]^-{k^\spc}\\
    \Alg Xt \ar[r]_-{\tresol} \ar[ur]|{\hole{\ds\weirdResol}\hole} & (\Alg Xt)^\spc
  }}
\]
We further define the natural transformation
\begin{equation}\label{eq:alpha}
  \vcenter{\xymatrix{
    \ar@{}[dr]|(0.7){\Uparrow\alpha} & A \ar[d]^-c\\
    A \ar@{=}[ur] \ar[r]_-{\resol} & A^\spc
  }}
\end{equation}
to be equal to
\[\xymatrix{
  A \ar@{=}[r] \ar[d]_-c \ar@{}[dr]|(0.7)\Uparrow & A \ar[d]^-c\\
  A^\aspc \ar[ur]|{\ev_{-1}} \ar[r]_-\res & A^\spc;
}\]
the transformation $\alpha^t\colon \tresol \Rightarrow c$ is defined similarly. Note that $k^\spc\alpha = \alpha^t k$.

\begin{definition}\label{def:cocomplete}
Let $g$ be a homotopy coherent comonad on an $\infty$-category $A$. An object $a\in A$ is \emph{$g$-cocomplete} if $\alpha_a$ is a colimiting cone. It is equivalent to ask that
\[\xymatrix{
  \ar@{}[dr]|(0.7){\Uparrow\alpha_a} & A \ar[d]^-c\\
  1 \ar[ur]^-a \ar[r]_-{\resol a} & A^\spc
}\]
be an absolute left lifting diagram, or that $\comma aA = \comma{\resol a}c$.
\end{definition}

\begin{definition}
If the colimit of $\resol a$ exists, we call it the \emph{$g$-cocompletion} of $a$ and notate it $\cocompl a$. It is characterized by $\comma{\resol a}{c} = \comma{\cocompl a}A$. Thus $a$ is $g$-cocomplete if and only if $\cocompl a$ exists and is equal to $a$.
\end{definition}

\begin{litcomp}
Our $g$-cocomplete objects correspond to the strongly $\bbK$-cocomplete objects of \cite[4.33]{Hess}.
\end{litcomp}

\begin{litcomp}
Suppose given model-categorical input. The \emph{derived $G$-cocompletion} $\cocompl[G] a$ of $a\in\cA$ is defined as the geometric realization of a Reedy cofibrant replacement of the simplicial object $\resol[G]R a$ given by the Blumberg-Riehl homotopical resolution. This is evidently compatible with our definition. The coaugmentation $\cocompl a \to a$ is modelled by the zig-zag
\[ \smash{\xymatrix{ \cocompl[G] a \ar[r] & Ra & \ar@{ >->}[l]_-\sim a}} \]
at the point-set level.
\end{litcomp}

We can now state the main theorem.

\begin{theorem}\label{thm:main}
The comparison functor $k$ is fully faithful on maps out of $a\in A$ if and only if $a$ is $g$-cocomplete. In particular, the restriction of $k$ to the full subcategory $A'$ of $g$-cocomplete objects is fully faithful, and $k$ is fully faithful if and only if $A'=A$.
\end{theorem}

\begin{corollary}
Suppose we are given model-categorical input. Then for each bifibrant $a\in\cA_\cf$, the comparison functor $A\morph^k \Alg Xt$ (in the world of $\infty$-categories) is fully faithful on maps out of $a\in A$ if and only if $a$ is weakly equivalent to its derived $G$-cocompletion $\cocompl[G] a$ (in the world of model categories).
\end{corollary}

\begin{remark}
We note that when the $\infty$-category $A$ has sufficiently many colimits, there is an easy proof of Theorem \ref{thm:main} not requiring any new results. Indeed, in this case the comparison functor has a left adjoint $\Alg Xt\morph^\magicAdj A$, given by $\magicAdj\defeq\mathord{\colim}\circ\weirdResol$, and the comonad $\ell k$ induced on $A$ is nothing but $g$-cocompletion $\cocompl\slot$. Thus \eqref{eq:alpha} defines an absolute left lifting diagram if and only if the counit of $\magicAdj\adjoint k$, which is the coagumentation from the cocompletion, is an isomorphism. But the counit of an adjunction is an isomorphism if and only if the right adjoint is fully faithful. Moreover, all these properties are determined pointwise \cite[6.1.8]{RV1}, so this argument proves the sharper condition as well. The results necessary to justify this line of reasoning in the $\infty$-categorical context are \cite[7.2.4]{RV2} and \cite[5.2.9]{RV5}.
\end{remark}

In general, we may not be able to define $\magicAdj$ on objects of $\Alg Xt$ not in the image of $k$. However, the ``non-representable left adjoint'' $\weirdResol$ turns out to suffice for the argument. In place of an adjunction $\magicAdj\adjoint k$, which would yield $\comma \magicAdj A = \comma{\Alg Xt}k$, we get

\begin{lemma}\label{lm:prt-adj}
We have $\comma \weirdResol c = \comma{\Alg Xt}k$.
\end{lemma}
\begin{proof}
Functors in each direction are given by
\[\xymatrix{
  \comma\weirdResol c \ar[r]^-{k^\spc} & \comma{\tresol}{ck} \\
  \comma\weirdResol\resol \ar[u]^-{\alpha} & \comma{\Alg Xt}k \ar[l]^-{\weirdResol} \ar@{=}[u]
}\]
where the identification comes from \cite[6.3.17]{RV2}. Since $\alpha\lambda$, $k^\spc$, and the identification all commute with the projections to $\Alg Xt$ and $A$, it follows from the 2-cell induction and 2-cell conservativity properties of commas (c.f.\ \cite[3.3.20]{RV1}) that these define inverse equivalences between $\comma\weirdResol c$ and $\comma{\Alg Xt}k$.
\end{proof}

\begin{proof}[Proof of Theorem \ref{thm:main}]
Let $a\in A$; there is a commutative diagram
\[\xymatrix@dr{
  \comma aA \ar[r]^-{\alpha} \ar[d]_-k & \comma{\resol a} c\\
  \comma{ka}k \ar@{=}[ur]
}\]
where the identification comes from Lemma \ref{lm:prt-adj}. It follows that $k$ is an equivalence if and only if $\alpha$ is so; but this is precisely the statement of the theorem.
\end{proof}

\subsection{Applications to descent}
\label{sub:descent}

We will apply the results from the previous section to the monadic formulation of descent, and to descent spectral sequences.

The previous section discussed monads. As is usual in category theory, we would like to obtain corresponding results for comonads ``by duality'', without having to repeat the arguments. The following construction achieves this in the $\infty$-cosmic setting.

\begin{definition}
Let $\cE$ be an $\infty$-cosmos. We define $\cE^\oop$ to be the simplicially enriched category with
\begin{itemize}
  \item the same objects as $\cE$, and
  \item mapping spaces given by $\map_{\cE^\oop}(A,B) = \map_\cE(A,B)^\op$.
\end{itemize}

Let $\cC$ be a 2-category. We define $\cC^\oop$ to be the 2-category with
\begin{itemize}
  \item the same objects as $\cC$, and
  \item mapping categories $\hom_{\cC^\oop}(A, B) = \hom_\cC(A,B)^\op$.
\end{itemize}
(The notation $\cC^{\mathrm{co}}$ is more standard, but has led to perverse uses of the term ``cofibration''.) Observe that $\htc{\cE^\oop} = \htc{\cE}^\oop$.
\end{definition}

We quickly summarize what this means for us. The ``walking comonad'' $\Cmd$ is the full subcategory of $\Adj$ on the object $-$. We now add $\pm$ subscripts to distinguish between our extension/restriction operations, writing them as
\[
  \adjnctn{\sSet^\Mnd}{\sSet^\Adj}{\lan_+}{\res_+}
  \quad\text{and}\quad
  \adjnctn{\sSet^\Cmd}{\sSet^\Adj}{\lan_-}{\res_-}.
\]
Given a homotopy coherent comonad $\ndo Ag$ coming from $\Cmd\morph^H\cE$ we have an $\infty$-category $\CoAlg Ag$ of homotopy coherent $g$-coalgebras, defined through weights by
\[\CoAlg Ag=\wLim{\res_-\Adj^+}H\Cmd\]
and producing a comonadic adjunction $\smash[b]{\adjnctn {\CoAlg Ag}A{u_g}{f_g}}$, with $u_gf_g=g$. If the comonad is induced from an adjunction $\ds\adjnctn XAfu$ associated to $\Adj\morph^T\cE$, then $\CoAlg Ag=\wLim{\lan_-\res_-\Adj^+}T\Adj$, and we have a cocomparison functor $X\morph^\kappa \CoAlg Ag$ fitting into a commutative diagram
\[\xymatrix{
  X \ar@(u,l)[]_t \ar[d]_-\kappa \ar@<1ex>[r]^-f \ar@{}[r]|-\bot
  &
  A \ar@(u,r)[]^g \ar@<1ex>[l]^-u \ar@{=}[d]
  \\
  \CoAlg Ag \ar@(d,l)^-{t_g} \ar@<1ex>[r]^-{u_g} \ar@{}[r]|-\bot
  &
  A\ar@<1ex>[l]^-{f_g}.
}\]

Let $\ndo Xt$ be a homotopy coherent monad in $\cE$ coming from a functor $\Mnd\morph^H\cE$. The monadic adjunction $\adjnctn X{\Alg Xt}{f^t}{u^t}$ is classified by $\wLim{\res_+\Adj^\slot}H\Mnd \colon\Adj\to\cE$, and induces a comonad $g^t=f^tu^t$ on $\Alg Xt$.

\begin{definition}
The \emph{$\infty$-category $\Desc Xt$ of descent data} for the monad $t$ is the $\infty$-category of $g^t$-coalgebras in $\Alg Xt$, $\Desc Xt \defeq \CoAlg{(\Alg Xt)}{g^t}$.
\end{definition}

Keeping in the spirit of the previous section, we would like a description of $\Desc Xt$ as a weighted limit. (We shan't need this, but the weight to use is not immediately obvious, and may be of use to posterity.) Let $\frW$ denote the subfunctor of $\Mnd^+=\res_+\Adj^+$ consisting of endomorphisms which factor through $-$ in $\Adj$; thus $\frW(+)=\Delta$. We can express this definition cleanly in terms of functor tensor products (which will be used in the proof) by $\frW = \comma\Cmd\Mnd \tnsr_\Cmd \comma+\Cmd$; we refer to \cite[\S\S4.1 and 4.3]{RiehlCHT} for an introduction to functor co/tensor products.

\begin{proposition}
$\Desc Xt = \wLim\frW H\Mnd$.
\end{proposition}
\begin{proof}
Expanding the definitions,
\[
  \Desc Xt \defeq \CoAlg{(\Alg Xt)}{g^t} = \wLim{\lan_-\res_-\Adj^+}{\wLim{\res_+\Adj^\slot}H\Mnd}\Adj.
\]
Applying the tensor hom-adjunction, this becomes
\[ \Desc Xt = \wLim{(\res_+\Adj^\slot) \tnsr_\Adj (\lan_-\res_-\Adj^+)}H\Mnd, \]
which is a description of $\Desc Xt$ as a single weighted limit over $\Mnd$. All that remains is to identify the weight; for this, we apply the tensor product formula for left Kan extension followed by the co-Yoneda lemma and get
\begin{align*}
  \Desc Xt &= \wLim{(\res_+\Adj^\slot) \tnsr_\Adj (\lan_-\res_-\Adj^+)}H\Mnd\\
  &= \wLim{ \comma\Adj\Mnd \tnsr_\Adj \comma\Cmd\Adj \tnsr_\Cmd \comma+\Cmd }H\Mnd\\
  &= \wLim{\comma\Cmd\Mnd \tnsr_\Cmd \comma+\Cmd}H\Mnd\\
  &= \wLim\frW H\Mnd. \qedhere
\end{align*}
\end{proof}
Thus the diagram of weights on the left induces the diagram of $\infty$-categories on the right.
\[
  \vcenter{\xymatrix@C=2em{
   \res_+\Adj^+ \ar@{}[r]|-\bot \ar@<-1ex>[r] & \res_+\Adj^- \ar@<-1ex>[l] \ar@{=}[d]\\
   \frW \ar@{ >->}[u] \ar@{}[r]|-\bot \ar@<-1ex>[r] & \res_+\Adj^- \ar@<-1ex>[l]
  }}
  \quad\rightsquigarrow\quad
  \vcenter{\xymatrix{
    X \ar[d]_-\delta \ar@{}[r]|{\bot} \ar@<1ex>[r]^{f^t} & \Alg Xt \ar@<1ex>[l]^{u^t} \ar@{=}[d]\\
    \Desc Xt \ar@{}[r]|{\bot} \ar@<1ex>[r]^{u_{g^t}} & \Alg Xt \ar@<1ex>[l]^{f_{g^t}}
  }}
\]

\begin{definition}
The $\infty$-monad $t$ satisfies \emph{descent} if $X\morph^\delta \Desc Xt$ is fully faithful. It satisfies \emph{effective descent} if $\delta$ is a weak equivalence.
\end{definition}

The dual of Theorem \ref{thm:main} immediately gives

\begin{proposition}
A $\infty$-monad $t$ on an $\infty$-category $X$ satisfies descent if and only if every object of $X$ is $t$-complete, if and only if every $x\in X$ is the totalization of its cosimplicial resolution $\coresol x$ determined by $t$.
\end{proposition}

\begin{corollary}
Suppose given model-categorical input. Then $t=N(T|_{\cX_\cf})$ satisfies descent (in the world of $\infty$-categories) if and only if every bifibrant $x\in\cX_\cf$ is weakly equivalent to its derived $T$-completion $\compl[T] x$ (in the world of model categories).
\end{corollary}

The remainder of this section is dual. Let $g$ be an $\infty$-comonad on the $\infty$-category $A$, coming from a functor $\Cmd\morph^H\cE$. The comonadic adjunction $\smash[b]{\adjnctn {\CoAlg Ag}A{u_g}{f_g}}$ induces a homotopy coherent monad $t_g=f_g u_g$ on $\CoAlg Ag$.

\begin{definition}
The \emph{$\infty$-category $\CoDesc Ag$ of codescent data} for the comonad $g$ is the $\infty$-category of $t_g$-algebras in $\CoAlg Ag$, $\CoDesc Ag \defeq \Alg{(\CoAlg Ag)}{t_g}$.
\end{definition}

Let $\frV$ denote the subfunctor of $\Cmd^- = \res_-\Adj^-$ consisting of endomorphisms which factor through $+$ in $\Adj$; thus $\frV=\comma\Mnd\Cmd\tnsr_\Mnd\comma-\Mnd$ and $\frV(-)=\Delta^\op$.

\begin{proposition}
$\CoDesc Ag = \wLim\frV H\Cmd$.
\end{proposition}

Thus the diagram of weights on the left induces the diagram of $\infty$-categories on the right.
\[
  \vcenter{\xymatrix@C=2em{
   \res_-\Adj^+ \ar@{=}[d] \ar@{}[r]|-\bot \ar@<-1ex>[r] & \res_-\Adj^- \ar@<-1ex>[l]\\
   \res_-\Adj^+ \ar@{}[r]|-\bot \ar@<-1ex>[r] & \frV \ar@{ >->}[u] \ar@<-1ex>[l]
  }}
  \quad\rightsquigarrow\quad
  \vcenter{\xymatrix{
    \CoAlg Ag \ar@{=}[d] \ar@{}[r]|{\bot} \ar@<1ex>[r]^{u_g} & A \ar@<1ex>[l]^{f_g} \ar[d]^-\gamma\\
    \CoAlg Ag \ar@{}[r]|{\bot} \ar@<1ex>[r]^{f^{t_g}} & \CoDesc Ag \ar@<1ex>[l]^{u^{t_g}}
  }}
\]

\begin{definition}
The $\infty$-comonad $g$ satisfies \emph{codescent} if $A\morph^\gamma \CoDesc Ag$ is fully faithful. It satisfies \emph{effective codescent} if $\gamma$ is a weak equivalence.
\end{definition}

\begin{proposition}
A $\infty$-comonad $g$ on an $\infty$-category $A$ satisfies codescent if and only if every object of $A$ is $g$-cocomplete, if and only if every $a\in A$ is the geometric realization of the simplicial resolution $\resol a$ given by $g$.
\end{proposition}

\begin{corollary}
Suppose given model-categorical input. Then $g=N(G|_{\cA_\cf})$ satisfies codescent (in the world of $\infty$-categories) if and only if every bifibrant $a\in\cA_\cf$ is weakly equivalent to its derived $G$-cocompletion $\cocompl[G] a$ (in the world of model categories).
\end{corollary}

\subsubsection{Spectral sequences}

Descent spectral sequences fall out easily in our setting. Our discussion follows \cite[\S5.3]{Hess} and is a trivial application of \cite[\S\S X.6--7]{BousfieldKan}. In particular, we refer to \cite[\S IX.5]{BousfieldKan} for treatment of convergence issues.

Let $t$ be an $\infty$-monad on an $\infty$-category $X$ with unit $\eta$, and assume that $X$ has all $t$-completions. Observe that a map $x\morph^\phi y$ gives rise to a cosimplicial pointed space $X(x, \coresol y)_{\coresol[\eta] \phi}$ whose totalization is $X(x, \compl y)_\phi$.

\begin{proposition}
A map $x\morph^\phi y$ in $X$ gives rise to a spectral sequence
\[ E^{r,s}_2 = \pi^r\pi_s X(x, \coresol y)_{\coresol[\eta] \phi} = \pi^r\pi_s\Desc Xt(\delta x, \delta \coresol y)_{\delta \circ(\coresol[\eta] \phi)} \]
which under suitable conditions converges to $\pi_* X(x, \compl y)_\phi$.
\end{proposition}

\begin{corollary}
Suppose given model-categorical input. For cofibrant $x$ and fibrant $y$ in $\cX$, a map $x\morph^\phi y$ gives rise to a spectral sequence
\[ E^{r,s}_2 = \pi^r\pi_s \Map_\cX(x, \coresol[T]Q y)_{\coresol[\zeta] \phi} \]
which under suitable conditions converges to $\pi_* \Map_\cX(x, \compl[T] y)_\phi$.
\end{corollary}

Dually, let $g$ be an $\infty$-comonad on an $\infty$-category $A$ with counit $\epsilon$, and assume that $A$ has all $g$-cocompletions. Observe that a map $a\morph^\psi b$ gives rise to a cosimplicial pointed space $A(\resol a, b)_{\psi\resol[\epsilon]}$ whose totalization is $A(\cocompl a, b)_\psi$.

\begin{proposition}
A map $a\morph^\psi b$ in $A$ gives rise to a spectral sequence
\[ E^{r,s}_2 = \pi^r\pi_s A(\resol a, b)_{\psi\resol[\epsilon]} = \pi^r\pi_s\CoDesc Ag(\gamma\resol a, \gamma b)_{\gamma \circ (\psi\resol[\epsilon])} \]
which under suitable conditions converges to $\pi_* A(\cocompl a, b)_\psi$.
\end{proposition}

\begin{corollary}
Suppose given model-categorical input. For cofibrant $a$ and fibrant $b$ in $\cA$, a map $a\morph^\psi b$ gives rise to a spectral sequence
\[ E^{r,s}_2 = \pi^r\pi_s \Map_\cA(\resol[G]R a, b)_{\psi\resol[\xi]} \]
which under suitable conditions converges to $\pi_* \Map_\cA(\cocompl[G] a, b)_\psi$.
\end{corollary}


\providecommand{\bysame}{\leavevmode\hbox to3em{\hrulefill}\thinspace}
\providecommand{\MR}{\relax\ifhmode\unskip\space\fi MR }
\providecommand{\MRhref}[2]{%
  \href{http://www.ams.org/mathscinet-getitem?mr=#1}{#2}
}
\providecommand{\href}[2]{#2}

\end{document}